\documentclass{svproc}
\usepackage{url}
\usepackage[colorlinks=true]{hyperref}
\hypersetup{urlcolor=blue, linkcolor=blue, citecolor=red, anchorcolor=blue]}
\usepackage{amsmath,amssymb,graphicx}
\usepackage{mathrsfs}
\newcommand{\R}{{\mathord{\mathbb R}}}
\renewcommand{\S}{\mathbb S}
\renewcommand{\(}{\left(}
\renewcommand{\)}{\right)}
\newcommand{\be}[1]{\begin{equation}\label{#1}}
\newcommand{\ee}{\end{equation}}
\newcommand{\ird}[1]{\int_{\R^d}{#1}\,dx}
\newcommand{\isd}[1]{\int_{\S^d}{#1}\,d\mu}
\newcommand{\irdg}[1]{\int_{\R^d}{#1}\,d\gamma}
\newcommand{\nrm}[2]{\left\|{#1}\right\|_{\mathrm L^{#2}(\R^d)}}
\newcommand{\nrms}[2]{\left\|#1\right\|_{\mathrm L^{#2}(\S^d)}}
\newcommand{\nrmng}[2]{\|#1\|_{\mathrm L^{#2}(\R^n,d\gamma)}}
\newcommand{\nrmG}[2]{\left\|#1\right\|_{\mathrm L^{#2}(\R^d,d\gamma)}}

\begin{document}

\mainmatter
\title{Stability results for Sobolev, logarithmic Sobolev, and related inequalities}
\titlerunning{Stability results for Sobolev type inequalities}
\author{Jean Dolbeault}
\authorrunning{J.~Dolbeault}
\tocauthor{Jean Dolbeault}

\institute{CEREMADE (CNRS UMR n$^\circ$ 7534),\\ PSL university, Universit\'e Paris-Dauphine,\\ Place de Lattre de Tassigny, 75775 Paris 16, France,\\
\email{dolbeaul@ceremade.dauphine.fr},\\ WWW home page:
\texttt{https://www.ceremade.dauphine.fr/\homedir dolbeaul}}
\maketitle
\begin{abstract}
Obtaining explicit stability estimates in classical functional inequalities like the Sobolev inequality has been an essentially open question for 30 years, after the celebrated but non-constructive result~\cite{MR1124290} of G.~Bianchi and H.~Egnell in 1991. Recently, new methods have emerged which provide some clues on these fascinating questions. The goal of the course is to give an introduction to the topic for some fundamental functional inequalities and present several methods that can be used to obtain explicit estimates.
\end{abstract}

The most classical Sobolev inequality on $\R^d$ with $d\ge3$ is 
\be{Sobolev}
\nrm{\nabla f}2^2\ge\mathsf S_d\,\nrm f{2^*}^2\quad\forall\,f\in\mathscr D^{1,2}(\R^d)
\ee
where $\mathscr D^{1,2}(\R^d)$ is the Beppo-Levi space, $\mathsf S_d$ denotes the sharp Sobolev constant and $2^*=2\,d/(d-2)$ is the critical Sobolev exponent. Equality in~\eqref{Sobolev} holds if and only if $f$ is in the $(d+2)$-dimensional manifold $\mathscr M$ of the \emph{Aubin--Talenti functions}
\[
g_{a,b,c}(x)=c\(a+|x-b|^2\)^{-\frac{d-2}2}\,,\quad a\in(0,\infty)\,,\quad b\in\R^d\,,\quad c\in\R\,.
\]
We refer to~\cite{frank2023sharp} for an up-to-date introduction to Ineq.~\eqref{Sobolev}. In 1991, G.~Bian\-chi and H.~Egnell proved in~\cite{MR1124290} that the \emph{deficit} in the Sobolev inequality, that is, the difference of the two sides of~\eqref{Sobolev}, satisfies the \emph{stability} inequality
\be{stab}
\nrm{\nabla f}2^2-\mathsf S_d\,\nrm f{2^*}^2\ge\kappa_d\,\min_{g\in\mathscr M}\nrm{\nabla(f-g)}2^2
\ee
for some constant~$\kappa_d>0$. The method of~\cite{MR1124290} uses compactness and provides no estimate on $\kappa_d$. Recently, several \emph{constructive} results have been obtained, \emph{i.e.}, results with explicit estimates of $\kappa_d$. The goal of these lecture notes is to make a short review based on two main strategies:
\begin{enumerate}
\item The strategy of G.~Bianchi and H.~Egnell made explicit in both of its steps: the global-to-local reduction and the local analysis, both relying on methods of the Calculus of Variations,
\item \emph{Entropy methods} and \emph{fast diffusion equations}, in which we can recognize also the same two steps, that are based on nonlinear parabolic flows and \emph{carr\'e du champ} techniques.
\end{enumerate}
Entropy methods usually result in better estimates of the stability constant, but to the price of some limitations. We will also state some results for related inequalities, like some families of Gagliardo-Nirenberg interpolation inequalities and logarithmic Sobolev inequalities, on the Euclidean space $\R^d$ or on the unit sphere $\S^d\subset\R^{d+1}$, and consider the use of the Legendre duality applied, for instance, to Hardy-Littlewood-Sobolev inequalities.

\section{Sobolev and Hardy-Littlewood-Sobolev inequalities}\label{Sec:Yamabe}

Let us start with a first example of fast diffusion flows applied to inequalities in the context of the Legendre duality. On $\R^d$ with $d\ge3$, the \emph{Hardy-Littlewood-Sobolev inequality}
\be{HLS}
\nrm g{\frac{2\,d}{d+2}}^2\ge\mathsf S_d\ird{g\,(-\Delta)^{-1}g}\quad\forall\,g\in\mathrm L^\frac{2\,d}{d+2}(\R^d)
\ee
is obtained by Legendre's duality from~\eqref{Sobolev}, as observed in~\cite{Lieb-83}. Indeed, by the Legendre transform, we obtain
\begin{align*}
&\sup_{f\in\mathscr D^{1,2}(\R^d)}\(\ird{f\,g}-\frac12\,\nrm f{2^*}^2\)=\frac12\,\nrm g{\frac{2\,d}{d+2}}^2\\
&\sup_{f\in\mathscr D^{1,2}(\R^d)}\(\ird{f\,g}-\frac12\,\nrm{\nabla f}2^2\)=\frac12\ird{g\,(-\Delta)^{-1}g}
\end{align*}
where $q=\frac{d+2}{d-2}$ is the H\"older conjugate of $2^*$, so that $\nrm{f^q}{\frac{2\,d}{d+2}}=\nrm f{2^*}^q$. It is not difficult to check that $\mathsf S_d$ is also the optimal constant in~\eqref{HLS}.

\begin{theorem}[\cite{MR3227280}]\label{Theorem:ExplicitGap0} Assume that $d\ge3$. There exists a constant $\mathscr C\in[\frac d{d+4},1)$ such that
\be{Sob-HLS}
\begin{aligned}
&\nrm{f^q}{\frac{2\,d}{d+2}}^2-\mathsf S_d\ird{f^q\,(-\Delta)^{-1}f^q}\\
&\hspace*{0.5cm}\le\mathscr C\,\mathsf S_d^{-1}\,\nrm f{2^*}^\frac 8{d-2}\(\nrm{\nabla f}2^2-\mathsf S_d\,\nrm f{2^*}^2\)\quad\forall\,f\in\mathscr D^{1,2}(\R^d)\,.
\end{aligned}
\ee
\end{theorem}
Theorem~\ref{Theorem:ExplicitGap0} can be considered as a stability result for~\eqref{Sobolev}, because the left-hand side in~\eqref{Sob-HLS}, which is the deficit associated with~\eqref{HLS}, measures a distance to $\mathscr M$. However, it is still a nonlinear functional and corresponds to a weaker norm than the usual norm in $\mathscr D^{1,2}(\R^d)$. With $\mathscr C=1$, there is an easy proof that goes as follows. Two integrations by parts show that
\[
\ird{|\nabla(-\Delta)^{-1}\,g|^2}=\ird{g\,(-\Delta)^{-1}\,g}
\]
and, if $g=f^q$ with $q=\frac{d+2}{d-2}$,
\[
\ird{\nabla f\cdot\nabla(-\Delta)^{-1}\,g}=\ird{f\,g}=\ird{g^q}=\ird{f^{2^*}}\,.
\]
Hence the expansion of the square
\[
0\le\ird{\left|\nrm f{2^*}^\frac4{d-2}\,\nabla f-\mathsf S_d\,\nabla(-\Delta)^{-1}\,g\right|^2}
\]
shows the result with $\mathscr C=1$. To prove that the inequality holds with $\mathscr C<1$ requires finer tools.

Let us consider the \emph{fast diffusion} equation
\be{Eqn:Yamabe}
\frac{\partial v}{\partial t}=\Delta v^m\,,\quad(t,x)\in\R^+\times\R^d
\ee
and compute the evolution of the deficit in~\eqref{HLS} defined as 
\[
-\,\mathsf H:=\nrm v{\frac{2\,d}{d+2}}^2-\mathsf S_d\,\ird{v\,(-\Delta)^{-1}v}\ge0\,.
\]
As in~\cite{MR2915466}, we observe that
\[
\frac 12\,\mathsf H'=\(\ird{v^\frac{2\,d}{d+2}}\)^\frac 2d\ird{\nabla v^m\cdot\nabla v^\frac{d-2}{d+2}}-\mathsf S_d\,\ird{v^{m+1}}
\]
where $v=v(t,\cdot)$ is a solution of~\eqref{Eqn:Yamabe}. With the choice $m=(d-2)/(d+2)$, that is, $m+1=2\,d/(d+2)=q$, which corresponds to the exponent in the \emph{Yamabe flow}, and $u=v^m$, we~have
\begin{multline*}
\frac 12\,\frac d{dt}\(\mathsf S_d\ird{v\,(-\Delta)^{-1}v}-\nrm v{\frac{2\,d}{d+2}}^2\)\\
=\(\ird{v^{m+1}}\)^\frac 2d\(\nrm{\nabla u}2^2-\mathsf S_d\,\nrm u{2^*}^2\)\ge0\,.
\end{multline*}
Eq.~\eqref{Eqn:Yamabe} admits a solution with \emph{separation of variables} vanishing at \hbox{$t=T>0$},
\[
\overline v_T(t,x)=c\,(T-t)^\alpha\,(g_\star(x))^\frac{d+2}{d-2}\,,\quad\alpha=\tfrac14\,(d+2)\quad\mbox{and}\quad c=\(4\,d\,\tfrac{d-2}{d+2}\)^\alpha\,,
\]
where $g_\star(x)=(1+|x|^2)^{-(d-2)/2}$ is the solution of $\Delta g_\star+d\,(d-2)\,g_\star^{(d+2)/(d-2)}=0$. Up to conformal transformations, this solution describes the asymptotic extinction profile of all solutions of~\eqref{Eqn:Yamabe} near the extinction time $T$.
\begin{theorem}[\cite{MR1857048}]\label{Thm:delPinoSaez}For any solution $v$ with initial datum $v_0\in\mathrm L^{2d/(d+2)}(\R^d)$, $v_0>0$, there exists $T>0$, $\lambda>0$ and $x_0\in\R^d$ such that
\[
\lim_{t\to T_-}(T-t)^{-\frac 1{1-m}}\,\sup_{x\in\R^d}(1+|x|^2)^{d+2}\,\left|\frac{v(t,x)}{\overline v(t,x)}-1\right|=0
\]
with $\overline v(t,x)=\lambda^{-(d+2)/2}\,\overline v_T\big(t,(x-x_0)/\lambda\big)$.\end{theorem}
The Hardy-Littlewood-Sobolev inequality~\eqref{HLS} amounts to $\mathsf H\le 0$ and appears as a consequence of the Sobolev inequality~\eqref{Sobolev} because $\mathsf H'\ge 0$ for any $t\in(0,T)$ and $\lim_{t\to T_-}\mathsf H(t)=0$ by Theorem~\ref{Thm:delPinoSaez}. By computing $\mathsf H''$, one can obtain refined estimates. Assume that $t\in[0,T)$ and let $\mathsf Y$ be such that
\[
\mathsf H(t)=-\,\mathsf Y(\mathsf J(t))\quad\mbox{with}\quad\mathsf J(t):=\ird{v(t,x)^\frac{2\,d}{d+2}}\,.
\]
With $\kappa_0=\mathsf H'_0/\mathsf J_0$, we obtain the differential inequality
\[
\mathsf Y'\(\mathsf S_d\,s^{1+\frac2d}+\mathsf Y\)\le\tfrac{d+2}{2\,d}\,\kappa_0\,\mathsf S_d^2\,s^{1+\frac4d}\,,\quad\mathsf Y(0)=0\,,\quad\mathsf Y(\mathsf J_0)=-\,\mathsf H_0\,.
\]
The relation $\mathsf Y'(\mathsf J)\,\mathsf '+\mathsf H'=0$ and elementary integrations by parts provide the following refinement of~\eqref{Sob-HLS} written with $\mathscr C=1$. 
\begin{proposition}[\cite{MR3227280}]\label{Thm:Improved} Assume that $d\ge3$. Then for any $f\in\mathscr D^{1,2}(\R^d)$ we have
\[
\begin{aligned}
&\nrm{f^q}{\frac{2\,d}{d+2}}^2-\mathsf S_d\ird{f^q\,(-\Delta)^{-1}f^q}\\
&\hspace*{0.5cm}\le\nrm f{2^*}^{2\,\frac{d+2}{d-2}}\,\varphi\(\mathsf S_d^{-1}\,\nrm f{2^*}^{-2}\(\nrm{\nabla f}2^2-\mathsf S_d\,\nrm f{2^*}^2\)\)
\end{aligned}
\]
with $q=(d+2)/(d-2)$ and $\varphi(s):=\sqrt{1+2\,s}-1$ for any $s\ge0$.
\end{proposition}
An asymptotic analysis of $v(t,\cdot)$ as $t\to T_-$ and some spectral estimates discard the case $\mathscr C=1$ in Theorem~\ref{Theorem:ExplicitGap0}, while the lower bound $\mathscr C=d/(d+4)$ is obtained by considering a sequence of test functions. See~\cite{MR3227280} for details. In dimension $d=2$ similar results holds in which~\eqref{Sobolev} and~\eqref{HLS} are respectively replaced by Onofri's inequality and the logarithmic Hardy-Littlewood-Sobolev inequality while~\eqref{Eqn:Yamabe} has to be replaced by the logarithmic fast diffusion equation using the results of~\cite{MR1357953,MR1606339,vazquez2023surveymassconservationrelated} instead of Theorem~\ref{Thm:delPinoSaez}. Notice that stability results based on duality are known for (logarithmic) Hardy-Littlewood-Sobolev inequalities from~\cite{MR3640893,carlen2024stabilitylogarithmichardylittlewoodsobolevinequality}.
The method based on~\eqref{Eqn:Yamabe} also applies to the fractional Sobolev inequality and the dual Hardy-Littlewood-Sobolev inequality as shown in~\cite{jankowiak2014fractionalsobolevhardylittlewoodsobolevinequalities}. Even if the optimal value of the stability constant $\mathscr C$ is not known, a nice feature of Theorem~\ref{Theorem:ExplicitGap0} is that we have a rather precise estimate of $\mathscr C$. However, a significant drawback of the method is that the deficit of the Hardy-Littlewood-Sobolev inequality corresponds to a weaker notion of distance to $\mathscr M$ than what can be expected.

\section{Stability, fast diffusion equation and entropy methods}

On $\R^d$, let us consider the \emph{fast diffusion equation} written in self-similar variables,
\be{RFD}
\frac{\partial v}{\partial t}+\nabla\cdot\Big(v\,\big(\nabla v^{m-1}-\,2\,x\big)\Big)=0\,.
\ee
Up to a $t$-dependent rescaling, it is equivalent to~\eqref{Eqn:Yamabe}. By \emph{fast diffusion}, we mean $m<1$ and reserve the expression \emph{porous medium} to the regime $m>1$ while the heat equation case corresponds to $m=1$ (in this case, Eq~\eqref{RFD} is known as the \emph{Fokker-Planck equation}, and $v^{m-1}$ should be replaced by $-\,\log v$). An introduction can be found in the two classical books~\cite{MR2176734,MR2282669}. Section~\ref{Sec:Yamabe} was devoted to $m=m_*:=(d-2)/(d+2)$ with $d\ge3$. Here we are interested in the range $m_c:=(d-2)/d<m<1$ in which nonnegative solution with initial data in $\mathrm L^1(\R^d)$ at $t=0$ exist for any $t\ge0$ according to~\cite{Herrero1985}, and more specifically to the range $m_1\le m<1$ where $m_1:=(d-1)/d$. Note that $m_*<m_c<m_1<1$ in dimension $d\ge2$. In any case,~\eqref{RFD} admits a stationary solution
\[
\mathcal B(x):=\(1+|x|^2\)^\frac1{m-1}
\]
which is usually called a \emph{Barenblatt solution} and has finite mass $M>0$ for any $m\in(m_c,1)$. If $m\in(m_c,1)$, the mass of any nonnegative solution $v$ is conserved and for simplicity, we shall assume that $\ird{v(t,x)}=M$ for any $t\ge0$. See~\cite{vazquez2023surveymassconservationrelated} for an interesting survey on mass conservation and the role of Barenblatt solutions in the range $(m_c,1)$. With this notation, we can introduce the \emph{generalized entropy} (or \emph{free energy}) and \emph{Fisher information} functionals respectively given~by
\begin{align*}
&\mathcal F[v]:=\frac1{m-1}\ird{\big(v^m-\mathcal B^m-m\,\mathcal B^{m-1}\,(v-\mathcal B)\big)}\,,\\
&\mathcal I[v]:=\frac m{1-m}\ird{v\left|\nabla v^{m-1}-\nabla\mathcal B^{m-1}\right|^2}\,.
\end{align*}
A straightforward computation shows that, for a solution of~\eqref{RFD}, we have
\be{ddtF}
\frac d{dt}\mathcal F[v(t,\cdot)]=-\,\mathcal I[v(t,\cdot)]\,.
\ee
\begin{theorem}[\cite{DelPino2002}]\label{Th:2002} Assume that $m\in(m_1,1)$ if $d=1$ or $d=2$, and $m\in[m_1,1)$ if $d\ge3$. If $v$ is a nonnegative solution of~\eqref{RFD} with initial datum $v_0\in\mathrm L^1(\R^d)$ of mass $M$ such that $\mathcal F[v_0]$ is finite, then
\be{FDE:rate}
\mathcal F[v(t,\cdot)]\le\mathcal F[v_0]\,e^{-\,4\,t}\quad\forall\,t\ge0\,.
\ee
\end{theorem}
The estimate~\eqref{FDE:rate} is equivalent to the \emph{entropy -- entropy production} inequality
\be{EEP}
\mathcal I[v]\ge 4\,\mathcal F[v]\,,
\ee
which can also be rewritten as a family of \emph{Gagliardo-Nirenberg-Sobolev inequalities}
\be{GNS}
\nrm{\nabla f}2^\theta\,\nrm f{p+1}^{1-\theta}\ge\mathscr C_{\mathrm{GNS}}(p)\,\nrm f{2p}
\ee
with optimal constant. The equivalence is obtained using
\[
p=\frac 1{2\,m-1}\quad\Longleftrightarrow\quad m=\frac{p+1}{2\,p}
\]
with $v=f^{2\,p}$ so that $v^m=f^{p+1}$ and $v\,|\nabla v^{m-1}|^2=(p-1)^2\,|\nabla f|^2$. The range $m\in(m_1,1)$ means that $2\,p\in(2,2^*)$ where $2^*=2\,d/(d-2)$ if $d\ge3$ and $2^*=+\infty$ otherwise. The case $m=m_1$, \emph{i.e.}, $2\,p=2^*$, is also covered if $d\ge3$.

Equality in~\eqref{GNS} is achieved by $f_\star(x)=\mathcal B(x)^\frac1{2p}=\(1+|x|^2\)^{-\frac1{p-1}}$, which also means that $4$ is the optimal constant in~\eqref{EEP} and it is also the best possible decay rate in~\eqref{FDE:rate}. The exponent $\theta=\frac{d\,(p-1)}{(d+2-p\,(d-2))\,p}$ is determined by the invariance under scalings and it is easy to check that $\theta=1$ for $p=2^*/2$ and $d\ge3$. If $d=2$, we obtain the Euclidean Onofri inequality in the limit as $p\to+\infty$. On the other hand, $\lim_{p\to1_+}\mathscr C_{\mathrm{GNS}}(p)=1$ so that~\eqref{GNS} becomes an equality for $p=1$. By taking the derivative with respect to $p$ at $p=1$, we obtain the scale invariant form of the scale invariant \emph{Euclidean logarithmic Sobolev inequality}
\be{LSI}
\frac d2\,\log\(\frac2{\pi\,d\,e}\ird{|\nabla f|^2}\)\ge\ird{|f|^2\,\log|f|^2}
\ee
for any function $f\in\mathrm H^1(\R^d)$ such that $\nrm f2=1$. Ineq.~\eqref{LSI} can be found in~\cite[Theorem~2]{MR479373},~\cite[Inequality~(2.3)]{MR0109101} and~\cite[Inequality~(26)]{MR1132315}. See \cite{1501} for further comments and~\cite{Gross75} for the classical Gaussian form and its Euclidean counterpart.

We obtain a \emph{linearized free energy} and a \emph{linearized Fisher information}
\[
\mathsf F[w]:=\frac m2\ird{w^2\,\mathcal B^{2-m}}\quad\mbox{and}\quad\mathsf I[w]:=m\,(1-m)\ird{|\nabla w|^2\,\mathcal B}
\]
by considering $f_\varepsilon:=\mathcal B\,(1+\varepsilon\,\mathcal B^{1-m}\,w)$ and taking the $O(\varepsilon^2)$ in, respectively $\mathcal F[f_\varepsilon]$ and $ \mathcal I[f_\varepsilon]$. The functionals $\mathsf F$ and $\mathsf I$ are related by a \emph{Hardy-Poincar\'e inequality}.
\begin{proposition}[\cite{Blanchet2009,Dolbeault2011a,BDNS2021}]\label{Prop:Hardy-Poincare} Let $m\in[m_1,1)$ if $d\ge3$, $m\in(1/2,1)$ if $d=2$, and $m\in(1/3,1)$ if $d=1$. If $w\in\mathrm L^2(\R^d,\mathcal B^{2-m}\,dx)$ is such that $\nabla w\in\mathrm L^2(\R^d,\mathcal B\,dx)$, $\ird{w\,\mathcal B^{2-m}}=0$, then
\be{Ineq:Hardy-Poincare}
\mathsf I[w]\ge4\,\alpha\,\mathsf F[w]
\ee
with $\alpha=1$. Additionally, under the center of mass condition
\be{Center:asymptotic}
\ird{x\,w\,\mathcal B^{2-m}}=0\,,
\ee
then the inequality holds with $\alpha=2-d\,(1-m)$.
\end{proposition}
The value of $\alpha$ in Proposition~\ref{Prop:Hardy-Poincare} is optimal and corresponds to the spectral gap of the linearized evolution operator $\mathcal L$ associated with~\eqref{RFD}, and to an improved spectral gap under Condition~\eqref{Center:asymptotic} in the second case. In the first case, we recover that $\inf\mathcal I/\mathcal F=\inf\mathsf I/\mathsf F=4$, which means that the worst decay rate is achieved in~\eqref{FDE:rate} in the \emph{asymptotic regime} as $t\to+\infty$. The eigenspace of $\mathcal L$ associated with $\alpha=1$ corresponds to the equality case in~\eqref{Ineq:Hardy-Poincare} and determines the spectral gap. See~\cite{Denzler2005} for a description of the spectrum of $\mathcal L$ and~\cite{Bonforte2010c} for an equivalence with the setting of Proposition~\ref{Prop:Hardy-Poincare}. Under Condition~\eqref{Center:asymptotic}, an improved entropy -- entropy production inequality is obtained, at least for functions which are close enough to the manifold of the Barenblatt functions. 
\begin{corollary}[\cite{Dolbeault2011a,Bonforte2010c,BDNS2021}]\label{Cor:Hardy-Poincare} Let $m\in(m_1,1)$ if $d\ge2$, $m\in(1/3,1)$ if $d=1$, $\eta=2\,(d\,m-d+1)$ and $\chi=m/(266+56\,m)$. If $\ird v=\mathrm M$, $\ird{x\,v}=0$ and
\be{B:voisinage}
(1-\varepsilon)\,\mathcal B\le v\le(1+\varepsilon)\,\mathcal B
\ee
for some $\varepsilon\in(0,\chi\,\eta)$, then
\[
\mathcal I[v]\ge(4+\eta)\,\mathcal F[v]\,.
\]
\end{corollary}
If $v$ solves~\eqref{RFD} with $m\in(m_1,1)$ and an initial datum $v_0$ with centred mass, \emph{i.e.}, such that \hbox{$\ird{x\,v_0(x)}=0$}, then $\ird{x\,v(t,x)}=0$ for any $t\ge0$ so that $w=(v-\mathcal B)\,\mathcal B^{m-2}$ satisfies~\eqref{Center:asymptotic}. If additionally $v_0$ satisfies~\eqref{B:voisinage}, then
\[
\mathcal F[v(t,\cdot)]\le\mathcal F[v_0]\,e^{-(4+\eta)\,t}\quad\forall\,t\ge 0\,.
\]
We have obtained an improved entropy~-- entropy production inequality and an improved decay rate of $\mathcal F[v(t,\cdot)]$ by comparing the nonlinear free energy $\mathcal F$ and the Fisher information $\mathcal I$ with their linearized counterparts $\mathsf F$ and $\mathsf I$, and by getting rid of the eigenspace corresponding to $\alpha=1$ in~\eqref{Ineq:Hardy-Poincare}. The result of Corollary~\ref{Cor:Hardy-Poincare} can be extended to $m\notin[m_1,1)$: see~\cite{Blanchet2009}. Under technical conditions, better rates can be achieved in some cases as in~\cite{Denzler2015}. One of the issues is that Condition~\eqref{B:voisinage} is rather restrictive and determines only a small neighbourhood of the set of the Barenblatt functions. Here we are going to take advantage of a property of \emph{relative uniform convergence} that goes back to~\cite{Bonforte2006}: Condition~\eqref{B:voisinage} is satisfied by any solution $v(t,\cdot)$ for $t\ge t_\star$ if the \emph{threshold time} $t_\star$ is taken large enough. Based on results of~\cite{Bonforte2019a,BDNS2021}, we have the following~result.
\begin{theorem}[\cite{BDNS2021,BDNS2023}]\label{Thm:RelativeUniform} Assume that $m\in(m_1,1)$ if $d\ge2$, $m\in(1/3,1)$ if $d=1$ and let $\varepsilon\in(0,1/2)$, small enough, and $A>0$ be given. There exists an explicit threshold time $t_\star\ge0$ such that, if $v$ is a solution of~\eqref{RFD} with nonnegative initial datum $v_0\in\mathrm L^1(\R^d)$ satisfying $\ird{v_0}=\ird{\mathcal B}=\mathrm M$ and
\be{A}
\mathsf A[v_0]=\sup_{r>0}r^\frac{d\,(m-m_c)}{(1-m)}\int_{|x|>r}v_0\,dx\le A<\infty\,,
\ee
and if $B(t,\cdot)$ solves~\eqref{Eqn:Yamabe} with initial datum $\mathcal B$, then
\[
\sup_{x\in\R^d}\left|\frac{v(t,x)}{\mathcal B(x)}-1\right|\le\varepsilon\quad\forall\,t\ge t_\star\,.
\]
\end{theorem}
Under Condition~\eqref{A}, we have an improved decay rate for any $t$ in the \emph{asymptotic time range} $[t_\star,+\infty)$.

The next question is to find an improvement in the \emph{initial time layer} $[0,t_\star)$ as well. A key ingredient is the \emph{carr\'e du champ} method for which we can refer to~\cite[Chapter~2]{BDNS2021}. We learn for instance that for some strictly convex function~$\psi$ with $\psi(0)=0$, $\psi'(0)=0$, there is an improved entropy~-- entropy production inequality
\[
\mathcal I-4\,\mathcal F\ge\psi(\mathcal F)\ge0\,,
\]
so that better decay rates can be expected if $\mathcal F$ is large, which is precisely what can be expected in the initial time layer. In fact, a key estimate of the \emph{carr\'e du champ} method is the simpler observation that
\[
\frac d{dt}\mathcal I[v(t,\cdot)]\le-\,4\,\mathcal I[v(t,\cdot)]
\]
if $v$ solves~\eqref{RFD}, which allows us to prove a \emph{backward in time} estimate. Let us define $\mathcal Q[v]:=\frac{\mathcal I[v]}{\mathcal F[v]}$ and notice that, using~\eqref{ddtF}, we have
\[
\frac d{dt}\mathcal Q[v(t,\cdot)]\le\mathcal Q[v(t,\cdot)]\,\big(\mathcal Q[v(t,\cdot)]-4\big)\,.
\]
\begin{lemma}[\cite{BDNS2021}]\label{prop.backward} Assume that $m>m_1$ and $v$ is a solution to~\eqref{RFD} with nonnegative initial datum $v_0$. If for some $\eta>0$ and $t_\star\ge0$, we have $\mathcal Q[v(t_\star,\cdot)]\ge4+\eta$, then
\[
\mathcal Q[v(t,\cdot)]\ge\,4+\zeta\quad\forall\,t\in[0,t_\star]\quad\mbox{with}\quad\zeta:=\frac{4\,\eta\,e^{-4\,t_\star}}{4+\eta-\eta\,e^{-4\,t_\star}}\,.
\]\end{lemma}

In the subcritical range with $m>m_1$, the first consequence of Lemma~\ref{prop.backward} is an improved decay rate of the free energy for the solution of~\eqref{RFD}.
\begin{corollary}[\cite{BDNS2021}]\label{Cor:Improvedrate} Let $m\in(m_1,1)$ if $d\ge2$, $m\in(1/2,1)$ if $d=1$, and $A>0$. If $v$ is a solution of~\eqref{RFD} with nonnegative initial datum $v_0\in\mathrm L^1(\R^d)$ such that $\ird{v_0}=\mathrm M$, $\ird{x\,v_0}=0$ and $v_0$ satisfies~\eqref{A}, then there is an explicit $\zeta>0$ such that
\[
\mathcal F[v(t,.)]\le\mathcal F[v_0]\,e^{-\,(4+\zeta)\,t}\quad\forall\,t\ge0\,.
\]
\end{corollary}
The constant $\zeta$ can be explicitly computed in terms of $d$, $m$ and $A$ using Corollary~\ref{Cor:Hardy-Poincare}, Theorem~\ref{Thm:RelativeUniform}, and Lemma~\ref{prop.backward}. The second consequence of Lemma~\ref{prop.backward} is a \emph{stability result for the entropy -- entropy production inequality}, which  reads as
\[
\mathcal I[v]-\,4\,\mathcal F[v]\ge\frac\zeta{4+\zeta}\,\mathcal I[v]
\]
where the left-hand side is the deficit in~\eqref{EEP}, while the distance in the right-hand side is measured by the Fisher information, which is a quantity that controls strong norms.

In the critical case $m=m_1$ corresponding to the Sobolev inequality for $d\ge3$, we also obtain a \emph{constructive stability result}. Let $\mathsf g(x):=(1+|x|^2)^{-(d-2)/2}$ be the Aubin--Talenti function such that $\mathsf g^{2^*}=\mathcal B$.
\begin{theorem}[\cite{BDNS2021}]\label{Thm:Main} Let $d\ge3$, $A>0$, and $p=2^*/2$. Then for any nonnegative $f\in\mathrm L^{p+1}(\R^d)$ such that $\nabla f\in\mathrm L^2(\R^d)$ and $|x|\,f^p\in\mathrm L^2(\R^d)$,
\[
\ird{\!(1,x, |x|^2)\,f^{2^*}}=\ird{\!(1,x, |x|^2)\,\mathsf g}\quad\mbox{and}\quad
\sup_{r>0}r^d\!\int_{|x|>r}f^{2^*}\,dx\le A\,,
\]
then we have
\[
\nrm{\nabla f}2^2-\mathsf S_d^2\,\nrm f{2^*}^2\ge\mathscr C\,\ird{\big|(p-1)\,\nabla f+f^p\,\nabla\mathsf g^{1-p}\big|^2}
\]
for some explicit $\mathscr C>0$ which depends only on $d$ and~$A$. 
\end{theorem}
Extending the subcritical result of Corollary~\ref{Cor:Improvedrate} to the critical case uses the introduction of a time delay, as shown in~\cite[Chapter~6]{BDNS2021}. One can also obtain a stability result in the subcritical range $m\in(m_1,1)$, that is, $1<p<2^*/2$ for Ineq.~\eqref{GNS}, but the statement is more delicate as several scales are involved. A remarkable feature in Theorem~\ref{Thm:Main} is that we measure the distance to a specific Aubin--Talenti function and not to the manifold $\mathscr M$.

Summing up, we have obtained a construction of an explicit stability result using the \emph{fast diffusion equation} and \emph{entropy methods} based on a parabolic regularity theory made explicit, spectral estimates and various features of the \emph{carr\'e du champ} method. However, estimates of the stability constants are way smaller than the expected values of the optimal stability constants and are available only under the restriction~\eqref{A} on the tails of the functions, which is inherent to the method, as one can prove that $t_\star$ diverges if $A\to+\infty$ in~\eqref{A}.

\section{Explicit stability results for Sobolev and log-Sobolev inequalities, with optimal dimensional dependence}\label{Sec:DEFFL}

In this section, we just state the main results and refer the reader to~\cite{https://doi.org/10.48550/arxiv.2209.08651} for detailed results and to~\cite{DEFFL2024} for a shorter introduction to the main tools. The most important fact is that the method gives an explicit value to the constant in the \emph{stability} inequality~\eqref{stab}, but again one can suspect that the lower bound is not very good, even if it has an optimal dimensional dependence.

\begin{theorem}[\cite{https://doi.org/10.48550/arxiv.2209.08651}]\label{DEFFL} There is a constant $\beta>0$ with an explicit lower estimate which does not depend on $d$ such that for all $d\ge 3$ and all $f\in \mathrm H^1(\R^d)$ we have
\be{BE}
\nrm{\nabla f}2^2-\mathsf S_d\,\nrm f{2^*}^2\ge\frac\beta d\,\inf_{g\in\mathscr M}\nrm{\nabla f-\nabla g}2^2\,.
\ee
\end{theorem}
Using the inverse stereographic projection we can rewrite the result on $\S^d$ as
\begin{multline*}
\nrms{\nabla F}2^2-\tfrac14\,d\,(d-2)\(\nrms F{2^*}^2-\nrms F2^2\)\\
\ge\frac\beta d\,\inf_{G\in\mathscr M(\S^d)}\(\nrms{\nabla F-\nabla G}2^2+\tfrac14\,d\,(d-2)\,\nrms{F-G}2^2\)\,.
\end{multline*}
Here we take on $\S^d$ the uniform probability measure. Let $v\in\mathrm H^1(\mathbb R^n,dx)$ with compact support, $d\ge n$ large enough and consider
\[
u_d(\omega)=v\Big(\tfrac1{r_d}\,\omega_1,\tfrac1{r_d}\,\omega_2,\ldots,\tfrac1{r_d}\,\omega_n\Big)\,,\quad r_d=\sqrt{\tfrac d{2\,\pi}}
\]
where $\omega=(\omega_1,\omega_2,\ldots,\omega_{d+1})\in\S^d\subset\R^{d+1}$. If $d\gamma=e^{-\pi\,|x|^2}\,dx$ denotes the Gaussian measure on $\R^n$, by taking the limit as $d\to+\infty$, we obtain
\begin{multline*}
\lim_{d\to+\infty}d\(\nrms{\nabla u_d}2^2-\tfrac14\,d\,(d-2)\(\nrms{u_d}p^2-\nrms{u_d}2^2\)\)\\
=\nrmng{\nabla u}2^2-\pi\int_{\R^n}|u|^2\,\log\(\frac{|u|^2}{\nrmng u2^2}\)d\gamma\,,
\end{multline*}
where the right-hand side is the deficit in the \emph{Gaussian logarithmic Sobolev inequality} of~\cite{Gross75}, which is equivalent to~\eqref{LSI}. We can also let $d\to+\infty$ in~\eqref{BE}. The delicate part is to bound the coefficients $(a,b,c)$ of the Aubin--Talenti functions which minimize the distance in~\eqref{stab} uniformly with respect to $d$. Because of the rescaling, the gradient term in the distance is lost and stability is measured only by a $\mathrm L^2(\R^n,d\gamma)$ norm. See~\cite{https://doi.org/10.48550/arxiv.2209.08651} for a detailed proof.
\begin{corollary}[\cite{https://doi.org/10.48550/arxiv.2209.08651}]\label{logsob} With $\beta>0$ as in Theorem~\ref{DEFFL}, we have
\begin{multline*}
\nrmng{\nabla u}2^2-\pi\int_{\R^n}|u|^2\,\log\(\frac{|u|^2}{\nrmng u2^2}\)d\gamma\\ \ge\frac{\beta\,\pi}2\,\inf_{a\in\R^d\!,\,c\in\R}\int_{\R^n}|u-c\,e^{a\cdot x}|^2\,d\gamma\,.
\end{multline*}
\end{corollary}
The proof of Theorem~\ref{DEFFL} relies on three main steps.
\begin{enumerate}
\item A sequence of nonnegative functions built using the competing symmetries of E.~Carlen and M.~Loss in~ \cite{CarlenLoss}, which alternate a Schwarz symmetric decreasing rearrangement and a transformation which corresponds, up to the stereographic projection, to a rotation on the sphere. This step reduces the stability problem to a local problem, on a neighbourhood of $\mathscr M$. We also use a continuous Steiner symmetrization in the discussion of one of the cases and an argument of~\cite{Christ}.
\item A Taylor expansion with explicitly controlled remainder terms. To obtain the $d$-dependence, one has to perform a delicate analysis which involves an argument of~\cite{MR908654}.
\item The result for sign-changing functions is deduced from the result for the positive and negative parts by a convexity estimate.
\end{enumerate}

For full details on the proof of Theorem~\ref{DEFFL}, on has to refer to~\cite{https://doi.org/10.48550/arxiv.2209.08651}. How much is lost in the first two steps is unclear, and this is why a proof based, \emph{e.g.}, on a fast diffusion equation would be a significant added value. There are several results which complement Theorem~\ref{DEFFL} and Corollary~\ref{logsob}. We know from~\cite{arXiv:2210.08482} that $\beta<4\,d/(d+4)$. The equality case in~\eqref{stab} written for the optimal value of $\kappa_d$ is achieved according to~\cite{arXiv:2211.14185}. A stability result for the Hardy-Littlewood-Sobolev inequality follows from the duality method of~\cite{MR3640893} and has been generalized in~\cite{chen2023stability} to cover various exponents. This result also proves an explicit stability estimate for the fractional Sobolev inequality, which uses duality but otherwise follows the very same lines of  proof as in~\cite{https://doi.org/10.48550/arxiv.2209.08651}.

The case of the logarithmic Hardy-Littlewood-Sobolev inequality is also obtained in~\cite{carlen2024stabilitylogarithmichardylittlewoodsobolevinequality} using duality with the Onofri inequality, whose explicit stability has been obtained in~\cite{MR3878729} by other methods.

\section{More explicit stability results for the logarithmic Sobolev and Gagliardo-Nirenberg inequalities on $\S^d$}\label{Sec:BDS}

Here the goal is to reconciliate the flow methods of Section 2 with stability estimates without restrictions as in Section 3 and get estimate which may have the right order of magnitude compared to the optimal stability constants. Let us consider the family of \emph{subcritical Gagliardo-Nirenberg inequalities}
\be{GNS:Sd}
\nrms{\nabla F}2^2\ge d\,\mathcal E_p[F]:=\frac d{p-2}\(\nrms Fp^2-\nrms F2^2\)\quad\forall\,F\in\mathrm H^1(\S^d,d\mu)
\ee
for any $p\in[1,2)\cup(2,2^*)$, with $2^*:=\frac{2\,d}{d-2}$ if $d\ge3$ and $2^*=+\infty$ if $d=1$ or $2$. If $d\ge3$, the limit case $p=2^*$ of~\eqref{GNS:Sd} is the Sobolev inequality on the sphere. As another limit case when $p\to2$, we obtain the \emph{logarithmic Sobolev inequality on the sphere}
\be{LSI:Sd}
\isd{|\nabla F|^2}\ge d\,\mathcal E_2[F]:=\frac d2\isd{F^2\,\log\(\frac{F^2}{\nrms F2^2}\)}\quad\forall\,F\in\mathrm H^1(\S^d,d\mu)
\ee
where $d\mu$ is the uniform probability on $\S^d$. These inequalities were obtained in \cite{Bakry1985,Beckner_1993,MR1134481}. While equality is always achieved by constant functions, the optimal constant $d$ can be identified with the first positive eigenvalue of the Laplace-Beltrami operator on $\S^d$. Let $X_1$ be the corresponding eigenspace of spherical harmonic functions and denote by $\mathrm \Pi_1$ the orthogonal projection onto $X_1$. If $\varphi\in X_1$ and $F_\varepsilon=1+\varepsilon\,\varphi$, a simple Taylor expansion shows that $\nrms{\nabla F_\varepsilon}2^2-d\,\mathcal E_p[F_\varepsilon]=O(\varepsilon^4)$, so that there is no hope to control the deficit from below by the square of a distance based on usual norm, at least in the direction of $X_1$.
\begin{theorem}[\cite{2210}]\label{Thm:sphere-GNS-stability} Let $d\ge1$ and $p\in(1,2^*)$. There is an explicit stability constant $\mathscr S_{d,p}>0$ such that, for any $F\in\mathrm H^1(\S^d,d\mu)$, we have
\begin{multline*}
\isd{|\nabla F|^2}-d\,\mathcal E_p[F]\\
\ge\mathscr S_{d,p}\(\frac{\nrms{\nabla\mathrm \Pi_1 F}2^4}{\nrms{\nabla F}2^2+\nrms F2^2}+\nrms{\nabla(\mathrm{Id}-\mathrm \Pi_1)\,F}2^2\)\,.
\end{multline*}
\end{theorem}
This results is a consequence of the stability analysis of~\cite{Frank_2022}, of explicit improvements of~\eqref{GNS:Sd} and~\eqref{LSI:Sd} obtained for instance in~\cite{DEKL2014,1504} based on~\cite{MR2381156,MR3229793}, which involve a generalized \emph{carr\'e du champ} method for the flow
\[
\frac{\partial u}{\partial t}=u^{-p\,(1-m)}\(\Delta u+(m\,p-1)\,\tfrac{|\nabla u|^2}u\)
\]
with exponents $m\in\big[m_-(d,p),m_+(d,p)\big]$ where
\[
m_\pm(d,p):=\frac1{(d+2)\,p}\(d\,p+2\pm\sqrt{d\,(p-1)\,\big(2\,d-(d-2)\,p\big)}\)\,,
\]
and of spectral properties taken from the proof of~\cite{MR1230930}. As in Section~\ref{Sec:DEFFL} and with the same notation, we can consider the large dimensional limit. 
\begin{theorem}[\cite{2302}]\label{Thm:GNStoPLS} Let $n\ge1$ and $u\in\mathrm H^1(\mathbb R^n,dx)$ with compact support. For any $p\in(1,2)$, we have
\begin{multline*}
\lim_{d\to+\infty}d\(\nrms{\nabla u_d}2^2-\tfrac d{2-p}\(\nrms{u_d}2^2-\nrms{u_d}p^2\)\)\\
=\nrmng{\nabla u}2^2-\tfrac1{2-p}\(\nrmng u2^2-\nrmng up^2\)\,.
\end{multline*}
\end{theorem}
In Theorem~\ref{Thm:sphere-GNS-stability} we have the satisfaction to obtain a stability result which distinguishes the $X_1$ eigenspace, is based on a fast diffusion equation and does not involve any restriction on the function space, but the estimate is not uniform as $p\to2^*$, as could be expected from Theorem~\ref{DEFFL}.

There is still much more to understand, particularly concerning~\eqref{LSI:Sd}. In Corollary~\ref{logsob}, the stability of the \emph{Gaussian logarithmic Sobolev inequality} is measured only by a $\mathrm L^2(\R^n,d\gamma)$ norm. It is known from~\cite{MR4305006,Indrei2024} that an $\mathrm H^1(\R^n,d\gamma)$ stability result cannot be expected without additional conditions. The quest for optimal conditions is open but as an illustration, we can conclude with a last result.
\begin{theorem}[\cite{2303}]\label{ThmMain2} Let $d\ge1$. For any $\varepsilon>0$, there is some explicit $\mathscr C>1$ depending only on $\varepsilon$ such that, for any $u\in\mathrm H^1(\R^d,d\gamma)$ with
\[
\irdg{\left(1,x\right)\,|u|^2}=(1,0)\,,\quad\irdg{|u|^2\,e^{\,\varepsilon\,|x|^2}}<\infty
\]
for some $\varepsilon>0$, then we have
\[
\nrmG{\nabla u}2^2\ge\frac{\mathscr C}2\irdg{|u|^2\,\log|u|^2}\,.
\]
Moreover, if $u$ is compactly supported in a ball of radius $R>0$, we know that the inequality holds with $\mathscr C=1+\frac{\mathscr C_\star(\mathsf K_\star)-1}{1+R^2\,\mathscr C_\star(\mathsf K_\star)}$, $\mathscr C_\star(\mathsf K)=1+\frac1{432\,\mathsf K}$ and $\mathsf K_\star:=\,\max\big(d,\frac{(d+1)\,R^2}{1+R^2}\big)$.
\end{theorem}
Here the result is stated in the form of an improved inequality under orthogonality constraints and we use the heat flow case corresponding to the limit $m=1$ of the fast diffusion flow. The proof is based on a stability result of~\cite{Fathi_2016} for a restricted class of functions, which contains the solution of the heat flow after an  initial time layer recently determined in~\cite{chen2021dimension}. For a compactly supported initial datum, we obtain a very explicit expression using a result of \cite{lee2003geometrical}. One can certainly object that $1/432\approx0.00231481$ does not provide us with a glorious improvement of the constant, but this is definitely a step in the understanding of \emph{explicit} stability results and, as a consequence, better estimates of decay rates in nonlinear evolution problems.

\bigskip\noindent{\bf Acknowledgment:} This work has been partially supported by the Project \emph{Conviviality} (ANR-23-CE40-0003) of the French National Research Agency.\\
\copyright\,\the\year\ by the author. This paper may be reproduced, in its entirety, for non-commercial purposes. \href{https://creativecommons.org/licenses/by/4.0/legalcode}{CC-BY 4.0}


\end{document}